\documentclass[12pt]{article}
\usepackage{amsmath}
\usepackage{amssymb}
\usepackage{amsmath,amssymb,amsbsy,amsfonts,amsthm,latexsym,
amsopn,amstext,amsxtra,euscript,amscd}

\begin{document}

\def\A{\mathbb{A}}
\def\B{\mathbf{B}}
\def \C{\mathbb{C}}
\def \F{\mathbb{F}}
\def \K{\mathbb{K}}

\def \Z{\mathbb{Z}}
\def \P{\mathbb{P}}
\def \R{\mathbb{R}}
\def \Q{\mathbb{Q}}
\def \N{\mathbb{N}}
\def \Z{\mathbb{Z}}

\def\B{\mathcal B}
\def\e{\varepsilon}

\def\cA{{\mathcal A}}
\def\cB{{\mathcal B}}
\def\cC{{\mathcal C}}
\def\cD{{\mathcal D}}
\def\cE{{\mathcal E}}
\def\cF{{\mathcal F}}
\def\cG{{\mathcal G}}
\def\cH{{\mathcal H}}
\def\cI{{\mathcal I}}
\def\cJ{{\mathcal J}}
\def\cK{{\mathcal K}}
\def\cL{{\mathcal L}}
\def\cM{{\mathcal M}}
\def\cN{{\mathcal N}}
\def\cO{{\mathcal O}}
\def\cP{{\mathcal P}}
\def\cQ{{\mathcal Q}}
\def\cR{{\mathcal R}}
\def\cS{{\mathcal S}}
\def\cT{{\mathcal T}}
\def\cU{{\mathcal U}}
\def\cV{{\mathcal V}}
\def\cW{{\mathcal W}}
\def\cX{{\mathcal X}}
\def\cY{{\mathcal Y}}
\def\cZ{{\mathcal Z}}

\def\f{\frac{|\A||B|}{|G|}}
\def\AB{|\A\cap B|}
\def \Fq{\F_q}
\def \Fqn{\F_{q^n}}

\def\({\left(}
\def\){\right)}
\def\fl#1{\left\lfloor#1\right\rfloor}
\def\rf#1{\left\lceil#1\right\rceil}
\def\Res{{\mathrm{Res}}}

\newcommand{\comm}[1]{\marginpar{
\vskip-\baselineskip \raggedright\footnotesize
\itshape\hrule\smallskip#1\par\smallskip\hrule}}

\newtheorem{lem}{Lemma}
\newtheorem{lemma}[lem]{Lemma}
\newtheorem{prop}{Proposition}
\newtheorem{proposition}[prop]{Proposition }
\newtheorem{thm}{Theorem}
\newtheorem{theorem}[thm]{Theorem}
\newtheorem{cor}{Corollary}
\newtheorem{corollary}[cor]{Corollary}
\newtheorem{prob}{Problem}
\newtheorem{problem}[prob]{Problem}
\newtheorem{ques}{Question}
\newtheorem{question}[ques]{Question}
\newtheorem{rem}{Remark}

\title{Sums of fractions modulo $p$}

\author{
{\sc C. A. D\'iaz} and {\sc M. Z. Garaev}}

\date{}

\maketitle

\begin{abstract}
Let $\F_p$ be the field of residue classes modulo a large prime $p$. The present paper is devoted to the problem of  representability of
elements of $\F_p$ as sums of fractions of the form $x/y$ with $x,y$ from short intervals of $\F_p$.
\end{abstract}



\section{Introduction}

Throughout the paper $\varepsilon$ is a small fixed positive constant, $p$ is a prime number sufficiently large in terms of $\varepsilon.$ As usual, $\F_p$ denotes the field of residue classes modulo $p$.
The elements of $\F_p$  we will frequently associate with the set $\{0,1,\ldots, p-1\}$. Given an integer $x$ coprime to $p$ (or an element $x$ from $\F_p^{*}=\F_p\setminus\{0\}$) we use $x^*$ or $x^{-1}$ to denote its multiplicative inverse modulo~$p$.

Let  $\lambda\in \F_p$ be fixed and let
$\cI$ and $\cJ$ be two intervals in $\F_p$. We assume that $\cI$ and $\cJ$ are nonzero, that is $\cI\not=\{0\},\, \cJ\not=\{0\}$. Motivated by the recent work of Shparlinski~\cite{Shp}, we consider the equation
\begin{equation}
\label{eqn:main}
\sum_{i=1}^{n}\frac{x_i}{y_i} =  \lambda,
\end{equation}
where $x_i, y_j$ are variables that run through the intervals $\cI$ and $\cJ$
respectively.  Using exponential sum estimates Shparlinski obtained an asymptotic formula for the number of solutions of general linear congruences. In the case of~\eqref{eqn:main} his results imply nontrivial estimates under
some conditions imposed on the cardinalities of $\cI$ and $\cJ$ (see Lemma~\ref{lem:Shp} below). In particular, if $n\ge 3$ and $|\cI|=|\cJ|>p^{n/(3n-2)+\varepsilon}$,
then the asymptotic formula obtained by Shparlinski becomes nontrivial for any fixed constant $\varepsilon>0$ (here and below,  for a given set $\cX$ we use $|\cX|$ to denote its cardinality).

In the present paper we consider the problem of solvability of~\eqref{eqn:main}. Our results are based on combinatorial and analytical tools.
Although we do not get an asymptotic formula for the number of solutions, our results give the solvability of~\eqref{eqn:main}
under  weaker conditions on the sizes of $|\cI|$ and $|\cJ|$.

\begin{theorem}
\label{thm:solv8} Let $\cI$ and $\cJ$ be intervals of $\F_p$ such
$$
|\cI|^2\cdot|\cJ| > p^{1+\varepsilon},\quad |\cI|\cdot|\cJ|^2 > p^{1+\varepsilon}.
$$
Then for any $\lambda\in\F_p$  the equation
\begin{equation}
\label{eqn:sum8}
\sum_{i=1}^{8}\frac{x_i}{y_i} =  \lambda
\end{equation}
has a solution with $(x_1,\ldots,x_{8})\in\cI^{8}$ and $(y_1,\ldots,y_{8})\in\cJ^{8}$.
\end{theorem}

From Theorem~\ref{thm:solv8} it follows, in particular, that for any $\varepsilon>0$ there is $\delta=\delta(\varepsilon)>0$ such that if $\cI$ and $\cJ$ are intervals of $\F_p$
with
$$
|\cI|>p^{1/3+\varepsilon},\quad |\cJ|>p^{1/3-\delta},
$$
then any element $\lambda\in\F_p$  can be represented in the form~\eqref{eqn:sum8} for some $(x_1,\ldots,x_{8})\in\cI^{8}$ and $(y_1,\ldots,y_{8})\in\cJ^{8}$.

\begin{theorem}
\label{thm:solv12} Let $\cI$ and $\cJ$ be nonzero intervals of $\F_p$ such that
$$
|\cJ|>p^{5/119},\quad |\cI|\cdot|\cJ|^{21/20}>p^{3/4+\varepsilon}.
$$
Then for any $\lambda\in\F_p$  the equation
\begin{equation}
\label{eqn:sum12}
\sum_{i=1}^{12}\frac{x_i}{y_i} =  \lambda
\end{equation}
has a solution with $(x_1,\ldots,x_{12})\in\cI^{12}$ and $(y_1,\ldots,y_{12})\in\cJ^{12}$.
\end{theorem}

From Theorem~\ref{thm:solv12} it follows, in particular, that for any $\varepsilon>0$ there is $\delta=\delta(\varepsilon)>0$ such that if $\cI$ and $\cJ$ are intervals of $\F_p$ with
$$
|\cI| >p^{9/40+\varepsilon}, \quad|\cJ| > p^{1/2-\delta},
$$
then any element $\lambda\in\F_p$  can be represented in the form~\eqref{eqn:sum12}
for some  $(x_1,\ldots,x_{12})\in\cI^{12}$ and $(y_1,\ldots,y_{12})\in\cJ^{12}$.

\begin{theorem}
\label{thm:solv4k} Let $k$ be a fixed positive integer constant, $\cI$ and $\cJ$ be intervals of $\F_p$ such that
$$
|\cI|\cdot|\cJ|^{2k/(k+1)}>p^{1+\varepsilon}.
$$
Then for any $\lambda\in\F_p$  the equation
$$
\sum_{i=1}^{4k}\frac{x_i}{y_i} =  \lambda
$$
has a solution with  $(x_1,\ldots,x_{4k})\in\cI^{4k}$ and $(y_1,\ldots,y_{4k})\in\cJ^{4k}$.
\end{theorem}

In particular, for any $\varepsilon>0$ there is $\delta=\delta(\varepsilon,k)>0$ such that if $\cI$ and $\cJ$ be intervals of $\F_p$ with
$$
|\cI|>p^{\frac{1}{k+1}+\varepsilon},\quad |\cJ|>p^{\frac{1}{2} -\delta},
$$
then  any element $\lambda\in\F_p$ can be representable in the form
$$
\sum_{i=1}^{4k}\frac{x_i}{y_i} =  \lambda,
$$
for some $(x_1,\ldots,x_{4k})\in\cI^{4k}$ and $(y_1,\ldots,y_{4k})\in\cJ^{4k}$.

It is to be mentioned that if the interval $\cJ$ starts from the origin and $|\cJ|>p^{\varepsilon}$, then there is a positive integer $n=n(\varepsilon)$ such for any
element $\lambda\in\F_p$  the equation
$$
\sum_{i=1}^{n}\frac{1}{y_i}=\lambda
$$
has a solution with $y_i\in\cJ$, see Shparlinski~\cite{Shp2}. However, the problem is still open for intervals $\cJ$ of arbitrary
positions.

\section{Lemmas}
Given sets $\cX\subset\F_p$ and $\cY\subset\F_p$, the product set $\cX\cY$ is defined by
$$
\cX\cY=\{xy;\, x\in \cX,\, y\in\cY\}.
$$
 For a positive integer $k$, the $k$-fold sum of $\cX$, is defined by
$$
k\cX=\{x_1+\ldots+x_k;\, x_i\in \cX\}.
$$
We also use the notation $\cX^{-1}=\{x^{-1};\, x\in \cX\setminus\{0\}\}.$

From the results of Glibichuk~\cite{Glib} it is known if $|\cX||\cY|>2p$ then $8\cX\cY=\F_p$. Here we need its version given by
Garaev and Garcia~\cite{GGarcia} (see also Garcia~\cite{Garcia} for even a more general statement).

\begin{lemma}
\label{lem:GGarcia} Let $\cA,\cB,\cC,\cD$ be subsets of
$\F_p^*$ such that
$$
|\cA||\cC|>(2+\sqrt 2)p,\quad |\cB||\cD|>(2+\sqrt 2)p.
$$
Then
$$
(2\cA)(2\cB)+(2\cC)(2\cD)=\F_p.
$$
\end{lemma}

We remark that the constant $2+\sqrt{2}$ that appears in the condition of the lemma can be substituted by a smaller one,  but we do not need it here.

Next, we need the following result from Cilleruelo and Garaev~\cite{CillGar} which is based on the idea of Heath-Brown~\cite{HB}.
\begin{lemma}
\label{lem:CillGar}
Let $\cJ$ be an interval in $\F_p$ and $\lambda\in \F_p^*$. Then the number $W_{\lambda}$ of solutions of the congruence
$$
xy=\lambda,\qquad x\in \cJ, \, y\in \cJ,
$$
satisfies
$$
W_{\lambda} < \frac{|\cJ|^{3/2+o(1)}}{p^{1/2}}+|\cJ|^{o(1)}.
$$
\end{lemma}

Observe that  for $\lambda\in\F_p^*$ the equation  $x^{-1}+y^{-1}=\lambda$ implies that
$$
(x-\lambda^{-1})(y-\lambda^{-1})=\lambda^{-2}.
$$
Hence, we have the following consequence of Lemma~\ref{lem:CillGar}.
\begin{corollary}
\label{cor:CillGar}
Let $\cJ$ be an interval in $\F_p$ and $\lambda\in \F_p^*$. Then the number $W_{\lambda}$ of the solutions of the congruence
$$
x^{-1}+y^{-1}=\lambda,\qquad x\in \cJ, \, y\in \cJ,
$$
satisfies
\begin{equation}
\label{eqn:CillGar}
W_{\lambda} < \frac{|\cJ|^{3/2+o(1)}}{p^{1/2}}+|\cJ|^{o(1)}.
\end{equation}
\end{corollary}
We recall that~\eqref{eqn:CillGar} is equivalent to the claim that for any $\varepsilon>0$ there exists $c=c(\varepsilon)>0$ such that
$$
W_{\lambda} <  c\Bigl(\frac{|\cJ|^{3/2+\varepsilon}}{p^{1/2}}+|\cJ|^{\varepsilon}\Bigr).
$$

We also need the following result of Bourgain and Garaev~\cite{BG}.

\begin{lemma}
\label{lem:BG}  Let $\cJ$ be an arbitrary nonzero interval in $\F_p$.  For any fixed positive integer constant $k$ the number $T_{k}$ of solutions of the congruence
\begin{equation}
\label{eqn:k=k}
y_1^{-1}+\ldots+y_k^{-1}=y_{k+1}^{-1}+\ldots+y_{2k}^{-1},\qquad y_1,\ldots,y_{2k}\in \cJ,
\end{equation}
satisfies
\begin{equation}
\label{eqn:thm kI*}
T_{k}<\Bigl(|\cJ|^{2k^2/(k+1)}+\frac{|\cJ|^{2k}}{p}\Bigr)|\cJ|^{o(1)}.
\end{equation}
\end{lemma}

Finally, we state the result of Shparlinski~\cite{Shp} which will be used to deal with
Theorem~\ref{thm:solv12} for relatively small intervals $\cJ$.
\begin{lemma}
\label{lem:Shp}  Let $\cI$ and $\cJ$ be two nonzero intervals in $\F_p$.  Then the number $R=R(\lambda,\cI,\cJ)$ of solutions of \eqref{eqn:main} with $x_i\in \cI$ and $y_i\in \cJ$
satisfies
$$
\Bigl|R-\frac{|\cI|^n|\cJ|^n}{p}\Bigr|<|\cI| |\cJ|\Bigl(|\cI|^{n-2}+(p|\cJ|)^{(n-2)/2}\Bigr)p^{o(1)}.
$$
\end{lemma}

\section{Proofs}

\subsection{Proof of Theorem~\ref{thm:solv8}}
We can assume that $|\cI|>10, |\cJ|>10$. Let $\cI_0\subset \F_p$ be an interval such that
$$
|\cI_0|>0.3|\cI|,\quad 2\cI_0=\cI_0+\cI_0\subset \cI.
$$
Such an interval obviously exists. Let $W_{\lambda}$ be the number of solutions of the congruence
$$
x^{-1}+y^{-1}=\lambda,\quad x\in\cJ, y\in\cJ.
$$
Using Corollary~\ref{cor:CillGar}, we have
$$
|\cJ|^2\ll \sum_{\lambda\in \cJ^{-1}+\cJ^{-1}} W_{\lambda}\le |\cJ^{-1}+\cJ^{-1}|\cdot \Bigl(\frac{|\cJ|^{3/2+o(1)}}{p^{1/2}}+|\cJ|^{o(1)}\Bigr).
$$
It follows that
$$
|\cJ^{-1}+\cJ^{-1}|>\min\Bigl\{|\cJ|^{2+o(1)}, \, p^{1/2}|\cJ|^{1/2+o(1)}\Bigr\}.
$$
Denote
$$
\cA=\cD=\cI_0\setminus\{0\},\quad \cB=\cC=\Bigl(\cJ^{-1}+\cJ^{-1}\Bigr)\setminus\{0\}.
$$
We have
\begin{eqnarray*}
|\cA||\cC|=|\cB||\cD| &\ge& 0.1|\cI_0|\cdot|\cJ^{-1}+\cJ^{-1}|\\  &\ge & \min\Bigl\{|\cI||\cJ|^{2+o(1)}, \, (p|\cI|^2|\cJ|)^{1/2+o(1)}\Bigr\}\\  &\ge &
p^{1+0.1\varepsilon}>4p.
\end{eqnarray*}
Thus, the condition of  Lemma~\ref{lem:GGarcia} is satisfied. Therefore, we get
$$
(2\cI_0)(4\cJ^{-1})+(2\cI_0)(4\cJ^{-1})=\F_p.
$$
Since $2\cI_0\subset \cI$, the result follows.

\subsection{Proof of Theorem~\ref{thm:solv12}} Let  $R$ be the number of solutions of the congruence~\eqref{eqn:sum12} with $x_i\in\cI,\,y_j\in\cJ$.
There are three cases to consider.

\smallskip

{\sc Case 1.} $p^{5/119}<|\cJ|<p^{15/37}$.

\smallskip

In view of Lemma~\ref{lem:Shp} applied with $n=12$, the number $R$ satisfies
$$
R>\frac{|\cI|^{12}|\cJ|^{12}}{p}-|\cI|^{11}|\cJ|p^{0.1\varepsilon} - |\cI||\cJ|^6  p^{5+0.1\varepsilon}.
$$
From the condition of the theorem it follows that
$$
|\cI|^{11}|\cJ| p^{0.1\varepsilon} < \frac{0.1|\cI|^{12}|\cJ|^{12}}{p},\quad |\cI||\cJ|^6 p^{5+0.1\varepsilon}< \frac{0.1|\cI|^{12}|\cJ|^{12}}{p}.
$$
Therefore, $R>0$ and the result follows in this case.
\smallskip

{\sc Case 2.} $|\cJ|>p^{5/8}$.

\smallskip

We fix a nonzero element $x_0\in \cI$ and denote by $R_1$ the number of solutions of the equation
$$
\sum_{i=1}^{12}y_i^{-1} = \lambda x_0^{-1},\quad y_i\in\cJ.
$$
It suffices to show that $R_1>0$. Let $\cJ_1=\cJ\setminus\{0\}.$ Expressing $R_1$ via exponential sums and following the standard procedure, we get
$$
\Bigl|R_1-\frac{|\cJ_1|^{12}}{p}\Bigr|\le \frac{1}{p}\sum_{a=1}^{p-1}\Bigl|\sum_{y\in \cJ_1}e_p(ay^*)\Bigr|^{12}.
$$
Here and below, we use the abbreviation $e_p(z)=e^{2\pi i z/p}$. By the well-known estimate for incomplete Kloosterman sums we have
$$
\max_{\gcd(a,p)=1}\Bigl|\sum_{y\in \cJ_1}e_p(ay^*)\Bigr|\ll 2p^{1/2}\log p.
$$
We also have
$$
\frac{1}{p}\sum_{a=0}^{p-1}\Bigl|\sum_{y\in \cJ_1}e_p(ay^*)\Bigr|^{2} = |\cJ_1|.
$$
Therefore,
$$
R_1>\frac{|\cJ_1|^{12}}{p}-2^{10}|\cJ_1|p^{5}(\log p)^{10}.
$$
Since $|\cJ_1|\ge |\cJ|-1>0.5p^{5/8}$, we get that $R_1>0$ and the result follows in this case.

\smallskip

{\sc Case 3.} $p^{15/37}<|\cJ|<p^{5/8}$.

\smallskip

Following the notation of Lemma~\ref{lem:BG}, we denote by $T_{k}$ the number of solutions of the congruence~\eqref{eqn:k=k}.
From the well-known application of the Cauchy-Schwarz inequality it follows that
\begin{equation}
\label{eqn:T3}
T_3\le (T_2T_4)^{1/2}.
\end{equation}
From Corollary~\ref{cor:CillGar} we easily obtain that
$$
T_2\le \Bigl(\frac{|\cJ|^{3/2+o(1)}}{p^{1/2}}+|\cJ|^{o(1)}\Bigr) |\cJ|^{2}.
$$
Since $|\cJ|>p^{15/37}>p^{1/3}$, we get that
\begin{equation}
\label{eqn:T2}
T_2\le\frac{|\cJ|^{7/2+o(1)}}{p^{1/2}}.
\end{equation}
Furthermore, by Lemma~\ref{lem:BG} and the condition $|\cJ|<p^{5/8}$, we have
$$
T_4<|\cJ|^{32/5 +o(1)}+\frac{|\cJ|^{8+o(1)}}{p}<|\cJ|^{32/5 +o(1)}.
$$
Combining this estimate with~\eqref{eqn:T3} and~\eqref{eqn:T2}, we obtain that
$$
T_3 < \frac{|\cJ|^{99/20+o(1)}}{p^{1/4}}.
$$
From the relationship between the number of solutions of a symmetric equation and the
cardinality of the corresponding set, we have
$$
|3\cJ^{-1}|=|\cJ^{-1}+\cJ^{-1}+\cJ^{-1}|\gg \frac{|\cJ|^6}{T_3},
$$
implying that
$$
|3\cJ^{-1}|\ge |\cJ|^{21/20}p^{1/4-0.1\varepsilon}.
$$
As in the proof of Theorem~\ref{thm:solv8}, let $\cI_0\subset \F_p$ be an interval such that
$$
|\cI_0|>0.3|\cI|,\quad 2\cI_0\subset \cI.
$$
Denote
$$
\cA=\cD=\cI_0\setminus\{0\},\quad \cB=\cC=\Bigl(3\cJ^{-1}\Bigr)\setminus\{0\}.
$$
We have
\begin{eqnarray*}
|\cA||\cC|=|\cB||\cD| &\ge& 0.1|\cI_0|\cdot|3\cJ^{-1}|\\  &\ge & |\cI||\cJ|^{21/20}p^{1/4-0.2\varepsilon}\\  &\ge &
p^{1+0.1\varepsilon}>4p.
\end{eqnarray*}
Thus, the condition of  Lemma~\ref{lem:GGarcia} is satisfied. Therefore, we get
$$
(2\cI_0)(6\cJ^{-1})+(2\cI_0)(6\cJ^{-1})=\F_p.
$$
Since $2\cI_0\subset \cI$, this concludes the proof of Theorem~\ref{thm:solv12}.

\subsection{Proof of Theorem~\ref{thm:solv4k}}

There are two cases to consider.

\smallskip

{\sc Case 1.} $|\cJ|>p^{(k+1)/2k}$.

\smallskip

We fix a nonzero element $x_0\in \cI$ and denote by $R_2$ the number of solutions of the equation
$$
\sum_{i=1}^{4k}y_i^{-1} = \lambda x_0^{-1},\quad y_i\in\cJ.
$$
It suffices to show that $R_2>0$. Denoting $\cJ_1=\cJ\setminus\{0\}$ and following exactly the same argument as in the Case 2 of Theorem~\ref{thm:solv12}, we get
$$
R_2 > \frac{|\cJ_1|^{4k}}{p}-2^{4k-2}|\cJ_1|p^{2k-1}(\log p)^{4k-2}.
$$
Since $|\cJ_1|\ge |\cJ|-1 > 0.5p^{(k+1)/2k}$, we have $R_2>0$ and the claim follows in this case.

\smallskip

{\sc Case 2.} $|\cJ|<p^{(k+1)/2k}$.

\smallskip

We recall that $T_k$ is the number of solutions of the congruence~\eqref{eqn:k=k}. From Lemma~\ref{lem:BG}
it follows that in our case we have the bound
$$
T_k<\Bigl(|\cJ|^{2k^2/(k+1)}+\frac{|\cJ|^{2k}}{p}\Bigr)|\cJ|^{o(1)}<|\cJ|^{2k^2/(k+1)+o(1)}.
$$
Hence, from the relationship between the number of solutions of a symmetric equation and the
cardinality of the corresponding set, we have
$$
|k\cJ^{-1}|> \frac{|\cJ|^{2k}}{|\cJ|^{2k^2/(k+1)+o(1)}}>|\cJ|^{2k/(k+1)}p^{-0.1\varepsilon}.
$$
Hence, denoting by  $\cI_0\subset \F_p$ an interval such that $|\cI_0|>0.3|\cI|,\quad 2\cI_0\subset \cI$ we verify that the the condition of Lemma~\ref{lem:GGarcia}
is satisfied with $\cA=\cC=\cI_0$ and $\cB=\cD=k\cJ^{-1}$. Thus, we get that
$$
\cI(2k\cJ^{-1})+\cI(2k\cJ^{-1})=\F_p
$$
which finishes the proof of Theorem~\ref{thm:solv4k}.

Address of the authors:

\vspace{0.3cm}

C. A. D\'iaz, Centro de Ciencias Matem\'{a}ticas, Universidad
Nacional Aut\'onoma de M\'{e}xico, C.P. 58089, Morelia,
Michoac\'{a}n, M\'{e}xico,

Email: {\tt cdiaz@matmor.unam.mx}

\vspace{0.3cm}

M.~Z.~Garaev, Centro de Ciencias Matem\'{a}ticas, Universidad
Nacional Aut\'onoma de M\'{e}xico, C.P. 58089, Morelia,
Michoac\'{a}n, M\'{e}xico.

Email: {\tt garaev@matmor.unam.mx}

\end{document}